\documentclass[11pt,reqno]{amsart}
\usepackage[a4paper,textwidth=150mm,textheight=225mm,
            centering,headsep=10mm,footskip=12mm]{geometry}
\usepackage[T1]{fontenc}
\usepackage[utf8]{inputenc}
\usepackage{lmodern}
\usepackage[final,protrusion=true,expansion=true]{microtype}
\usepackage{amssymb,mathtools,mathrsfs}
\usepackage{xspace}
\usepackage{xcolor}
\definecolor{linkblue}{RGB}{20,50,120}
\definecolor{citegreen}{RGB}{20,100,60}
\usepackage[colorlinks=true,linkcolor=linkblue,citecolor=citegreen,
            urlcolor=linkblue,bookmarksnumbered=true]{hyperref}
\hypersetup{
  pdftitle={Characterisations of symmetric Birkhoff-James orthogonality via numerical ranges},
  pdfauthor={Janko Bra\v{c}i\v{c}, Cristina Diogo},
  pdfsubject={Spatial numerical range},
}
\theoremstyle{plain}
\newtheorem{theorem}{Theorem}[section]
\newtheorem{proposition}[theorem]{Proposition}
\newtheorem{corollary}[theorem]{Corollary}
\newtheorem{lemma}[theorem]{Lemma}
\theoremstyle{definition}
\newtheorem{remark}[theorem]{Remark}
\newtheorem{example}[theorem]{Example}
\newtheorem{definition}[theorem]{Definition}
\renewcommand{\theenumi}{\roman{enumi}} 
\renewcommand{\labelenumi}{(\theenumi)}
\numberwithin{equation}{section}
\newcommand{\pA}{\mathcal{A}}
\newcommand{\pB}{\mathcal{B}}
\newcommand{\pD}{\mathcal{D}}
\newcommand{\eM}{\mathscr{M}}
\newcommand{\eN}{\mathscr{N}}
\newcommand{\eH}{\mathscr{H}}
\newcommand{\eS}{\mathscr{S}}
\newcommand{\eX}{\mathscr{X}}
\newcommand{\eY}{\mathscr{Y}}
\newcommand{\eZ}{\mathscr{Z}}
\newcommand{\bC}{\mathbb{C}}
\newcommand{\bD}{\mathbb{D}}
\newcommand{\bR}{\mathbb{R}}
\newcommand{\re}{\operatorname{Re}}
\newcommand{\im}{\operatorname{Im}}
\newcommand{\dist}{\operatorname{dist}}
\newcommand{\conv}{\operatorname{conv}}
\newcommand{\SH}{\eS_{\!\!\eH}}
\newcommand{\SX}{\eS_{\!\!\eX}}
\newcommand{\PiX}{\Pi_{\!\eX}}
\newcommand{\SXd}{\eS_{\!\!\eX^*}}
\newcommand{\soi}{\textup{S0I}\xspace}
\newcommand{\bj}{\perp_{_{\mathrm{BJ}}}}
\newcommand{\tr}{\operatorname{tr}}
\newcommand{\Wr}{W^{\mathrm{r}}}
\begin{document}
%
%
\title[Characterisations of symmetric Birkhoff--James orthogonality]{Characterisations of symmetric Birkhoff--James orthogonality via 
numerical ranges}
\author[J. Bra\v{c}i\v{c}]{Janko Bra\v{c}i\v{c}}
\address{Faculty of Natural Sciences and Engineering, University of Ljubljana,
A\v{s}ker\v{c}eva c. 12, SI-1000 Ljubljana, Slovenia}
\email{janko.bracic@ntf.uni-lj.si}

\author[C. Diogo]{Cristina Diogo}
\address{ISCTE - Instituto Universitário de Lisboa, Avenida das For\c cas Armadas 1649-026 Lisboa, Portugal,
and Centro de Análise Matemática, Geometria e Sistemas Dinâmicos, Departamento de Matemática, Instituto Superior Técnico, Universidade de Lisboa, Av. Rovisco Pais 1, 1049-001 Lisbon, Portugal}
\email{cristina.diogo@iscte-iul.pt}

%
%
\begin{abstract}
Birkhoff--James orthogonality in a normed space is in general not symmetric. For a complex Banach space $\eX$ of dimension
at least $2$, we prove that the following assertions are equivalent: Birkhoff--James orthogonality in $\eX$ is symmetric; the
spatial numerical range of every operator on $\eX$ coincides with its reverse numerical range, defined by reversing the
orthogonality; $\eX$ has the simultaneous zero inclusion property, that is, $0$ lies in the closure of
the numerical range of an invertible operator exactly when it lies in that of its inverse;
$\pB(\eX)$ has the analogous property for the algebra numerical range; and a simple inequality for supporting functionals
holds. This confirms a conjecture of Ghosh. We also introduce the asymmetry constant $R(\eX)\in[1,2]$, as the supremum of the
norms of the projections $I-x\otimes\xi$, where $x$ is a unit vector and $\xi$ is a supporting functional at $x$. We show
that $R(\eX)^{-1}$ is the radius of the largest ball around the identity on which the simultaneous zero inclusion property
holds, and that every number in $[1,2]$ is the asymmetry constant of some two-dimensional $\ell^p$ space.
\end{abstract}
\subjclass[2020]{46B20, 47A12, 46C15}
\keywords{Birkhoff--James orthogonality; symmetry of orthogonality; characterisation of inner product spaces;
spatial numerical range; reverse numerical range; simultaneous zero inclusion property; asymmetry constant}
%
\maketitle
\section{Introduction} \label{sec01}

A vector $x$ in a normed space $\eX$ is Birkhoff--James orthogonal to $y$, written $x\bj y$, if $\|x+\lambda y\|\geq\|x\|$
for all scalars $\lambda$ \cite{Bir,Jam47b}. In an inner product space this is the usual orthogonality, but in general the
relation is not symmetric. By the theorem of James and Day \cite{Jam47a,Day}, which holds also for complex spaces
\cite{Boh,AGKRZ}, a space of dimension at least $3$ is an inner product space if and only if Birkhoff--James orthogonality
in it is symmetric. In dimension two there are other examples: the Radon planes in the real case (see \cite{MS,AMW}) and the
complex planes constructed by Oman \cite{Oma}.

In \cite{BD} we introduced the simultaneous zero inclusion (\soi) property of a normed space $\eX$: for every invertible
operator $T$ on $\eX$, the number $0$ lies in the closure of the spatial numerical range $W(T)$ \cite{Lum} if and only if
it lies in the closure of $W(T^{-1})$. Poon \cite{Poo} observed that the symmetry of Birkhoff--James orthogonality implies the
\soi property, and he asked whether the converse holds. Ghosh \cite{Gho} proved the converse for several classes of spaces
and conjectured that it holds in general.

Our main result, Theorem \ref{theo34}, confirms this conjecture. For a complex Banach space $\eX$ of dimension at least $2$ it
shows that the following assertions are equivalent: Birkhoff--James orthogonality in $\eX$ is symmetric; $W(S)=\Wr(S)$ for
every operator $S$, where the reverse numerical range $\Wr(S)$ is obtained by reversing the orientation of the orthogonality
in the definition of $W(S)$; $\eX$ has the \soi property; the algebra $\pB(\eX)$ has the analogous property for the algebra numerical 
range; and $\big|1-\langle x,\eta\rangle\langle u,\xi\rangle\big|\leq1$ for all unit vectors $x,u$ and all supporting functionals $\xi$ at $x$ 
and $\eta$ at $u$. The last condition is a complex counterpart of a characterisation of real Radon planes due to James \cite{Jam47b} and 
Martini and Swanepoel \cite{MS}. The proof rests on two observations. The first is that, if the last condition fails, a rank-one perturbation 
$T=I-\nu\,u\otimes\eta$ of the identity has its numerical range in a half-plane which does not contain $0$, while $0\in W(T^{-1})$. 
The second observation is that, if $u\bj y$ but $y\not\bj u$, a small perturbation of $y$ in the direction of $u$, together with a 
weak$^*$ compactness argument, produces pairs which violate the condition.

In Section \ref{sec04} we measure the asymmetry of Birkhoff--James orthogonality. The asymmetry constant $R(\eX)$ is the supremum 
of the norms of the projections $I-x\otimes\xi$, where $x$ is a unit vector and $\xi$ is a supporting functional at $x$. The constant lies 
in $[1,2]$, and it is equal to $1$ if and only if Birkhoff--James orthogonality is symmetric. Its reciprocal is an infimum of minimum gaps 
in the sense of Kato \cite{Kat}, and it is analogous to a constant of Ji and Wu \cite{JW}. We show in Theorem \ref{theo46} that
$R(\eX)^{-1}$ is the largest radius of a ball around $I$ on which the \soi property holds locally, and that $R(\eX)$ is the
best constant in the inclusion $\Wr(S)\subseteq\overline{D(0,R(\eX)\|S\|)}$. We have $R(\eX)=2$ for $\ell^1$, $\ell^\infty$, $C(K)$ and 
many other classical spaces. We show that every number in $[1,2]$ is the asymmetry constant of some space $\ell^p(2)$.

\section{Preliminaries} \label{sec02}

Throughout this paper, let $\eX$ be a complex Banach space of dimension at least $2$. By $\pB(\eX)$ we denote the
Banach algebra of all bounded linear operators on $\eX$. The dual space is $\eX^*$, and the pairing between $\eX$ and
$\eX^*$ is given by $\langle x,\xi\rangle=\xi(x)$. Let $\SX$ and $\SXd$ be the unit spheres of $\eX$ and $\eX^*$,
respectively.

\subsection*{Spatial numerical range and the \soi property}
For $x\in\SX$, let
$$J(x)=\{\xi\in\eX^*;\;\|\xi\|\leq1\;\text{and}\;\langle x,\xi\rangle=1\}.$$
By the Hahn--Banach theorem, $J(x)$ is non-empty. If $\xi\in J(x)$, then $1=|\langle x,\xi\rangle|\leq\|x\|\|\xi\|=\|\xi\|$,
and therefore $\|\xi\|=1$, that is, $J(x)\subseteq\SXd$. The elements of $J(x)$ are called the {\em supporting functionals} at $x$. 
The set $J(x)$ is the intersection of the closed unit ball of $\eX^*$ with the hyperplane $\{\xi\in\eX^*;\;\langle x,\xi\rangle=1\}$. 
Both sets are convex. The closed unit ball is weak$^*$ compact by the Banach--Alaoglu theorem. The hyperplane is weak$^*$ closed, 
since the mapping $\xi\mapsto\langle x,\xi\rangle$ is weak$^*$ continuous. Therefore $J(x)$ is a weak$^*$ compact convex set.

Let $\PiX=\{(x,\xi);\;x\in\SX,\;\xi\in J(x)\}$. The {\em spatial numerical range} of $T\in\pB(\eX)$ (see \cite{Lum}) is
$$W(T)=\{\langle Tx,\xi\rangle;\;(x,\xi)\in\PiX\}.$$
Since $|\langle Tx,\xi\rangle|\leq\|T\|$ for $(x,\xi)\in\PiX$, $W(T)$ is a non-empty subset of the closed disc
$\overline{D(0,\|T\|)}$. Moreover, the spectrum of $T$ is contained in the closure of the spatial numerical range, that is, 
$\sigma(T)\subseteq\overline{W(T)}$. We refer the reader to \cite{BD1} for more details about numerical ranges.

In \cite[Definition 1.1]{BD} we introduced the \soi property for finite-dimensional normed spaces. We extend this
definition as follows.

\begin{definition} \label{def21}
The {\em \soi property holds locally} at an invertible operator $T\in\pB(\eX)$ if
$0\in\overline{W(T)}\iff0\in\overline{W(T^{-1})}$. If the \soi property holds locally at every invertible operator,
then {\em $\eX$ has the \soi property}.
\end{definition}

If $\eX$ is finite-dimensional, then $\PiX$ is compact, and $W(T)$ is closed. Hence Definition \ref{def21} agrees
with \cite[Definition 1.1]{BD} in this case. Poon \cite[Definition 1.4]{Poo} and Ghosh \cite{Gho} define the \soi
property without closures, also for infinite-dimensional spaces. By Theorem \ref{theo34} and Remark \ref{rem37}~(i), both
versions are equivalent to the symmetry of Birkhoff--James orthogonality. Hence $\eX$ has the \soi property in the sense of 
Definition \ref{def21} if and only if it has the \soi property in the sense of Poon and Ghosh.

\subsection*{Birkhoff--James orthogonality}
A vector $x\in\eX$ is {\em Birkhoff--James orthogonal} to $y\in\eX$, written $x\bj y$, if $\|x+\lambda y\|\geq\|x\|$ for all 
$\lambda\in\bC$ (see \cite{Bir} and \cite[Definition 1.2]{Jam47b}). The relation is {\em symmetric} if $x\bj y$ implies $y\bj x$ for all 
$x,y\in\eX$. The following characterisation was proved by James \cite[Theorem 2.1]{Jam47b} for real spaces. The same proof
works for complex spaces (see \cite[Lemma 1.5]{Poo}).

\begin{lemma} \label{lem22}
Let $x\in\SX$ and $y\in\eX$. Then $x\bj y$ if and only if there is $\xi\in J(x)$ with $\langle y,\xi\rangle=0$.\qed
\end{lemma}

It follows from the lemma that $0\in W(T)$ if and only if $x\bj Tx$ for some non-zero $x\in\eX$, since Birkhoff--James orthogonality is 
homogeneous. This observation leads to the following result of Poon \cite[Proposition 1.6 and Corollary 1.7]{Poo}.

\begin{lemma} \label{lem23}
Suppose that Birkhoff--James orthogonality is symmetric in $\eX$. If $T\in\pB(\eX)$ is invertible, then
$0\in W(T)\iff 0\in W(T^{-1})$.\qed
\end{lemma}

\subsection*{The algebra numerical range and the algebraic \soi property}
We always assume that $\pA$ is a complex Banach algebra with unit $e$ such that $\|e\|=1$. Let
$\pD(\pA)=\{\varphi\in\pA^*;\;\varphi(e)=1=\|\varphi\|\}$ be its set of states. The {\em algebra numerical range} of $a\in\pA$ is
$$V(a)=\{\varphi(a);\;\varphi\in\pD(\pA)\}.$$
It is a compact convex set containing $\sigma(a)$ (see \cite{BD1} and \cite[Theorem 1]{SW}). For an operator $T\in\pB(\eX)$
one has $V(T)=\overline{\conv(W(T))}$ \cite[Theorem 6]{SW}. Since $V(a)$ is compact and convex, $0\notin V(a)$ if
and only if there exist $\theta\in[0,2\pi)$ and $\varepsilon>0$ such that
\begin{equation} \label{eq21}
\re\big(e^{i\theta}\lambda\big)\geq\varepsilon\qquad\text{for all}\;\lambda\in V(a).
\end{equation}
For $T\in\pB(\eX)$ this means that $\overline{W(T)}$ lies in a closed half-plane which does not contain $0$. We
say that such an operator satisfies the {\em half-plane condition}. By a result of Williams \cite[Theorem 1]{Wil},
$0\in V(a)$ if and only if $\|a-\lambda e\|\geq|\lambda|$ for all $\lambda\in\bC$, that is, if and only if $e\bj a$.

\begin{definition} \label{def24}
A unital Banach algebra $\pA$ has the {\em algebraic \soi property} if
$$0\in V(a)\iff0\in V(a^{-1})\qquad\text{for every invertible}\;a\in\pA.$$
\end{definition}

Every Hilbert space has the \soi property (see the paragraph before Definition 1.1 in \cite{BD}). The following quantitative
version is a consequence of \cite[Theorem 3.6]{CL}, where $\|T\|$ is even replaced by the numerical radius of $T$. We include a
short proof, which also covers the algebraic \soi property.

\begin{lemma} \label{lem25}
Let $\eH$ be a complex Hilbert space. If $T\in\pB(\eH)$ is invertible, then
$$\dist\big(0,W(T^{-1})\big)\geq\frac{\dist\big(0,W(T)\big)}{\|T\|^2}\qquad\text{and}\qquad
\dist\big(0,W(T)\big)\geq\frac{\dist\big(0,W(T^{-1})\big)}{\|T^{-1}\|^2}.$$
In particular, $\eH$ has the \soi property, and $\pB(\eH)$ has the algebraic \soi property.
\end{lemma}

\begin{proof}
Every $y\in\SH$ can be written as $y=\frac{1}{\|Tx\|}Tx$ with $x=\frac{1}{\|T^{-1}y\|}T^{-1}y\in\SH$. Note that
$$\langle T^{-1}y,y\rangle=\frac{1}{\|Tx\|^2}\,\overline{\langle Tx,x\rangle}.$$
Since $\|Tx\|\leq\|T\|$, we get $|\langle T^{-1}y,y\rangle|\geq\|T\|^{-2}\dist(0,W(T))$. This is the first
inequality. The second one follows by interchanging $T$ and $T^{-1}$. Hence $0\in\overline{W(T)}$ if and only if
$0\in\overline{W(T^{-1})}$. By the Toeplitz--Hausdorff theorem, $W(T)$ is convex. Hence $V(T)=\overline{W(T)}$ by
\cite[Theorem 6]{SW}, and the same holds for $T^{-1}$. This gives the algebraic \soi property of $\pB(\eH)$.
\end{proof}

\subsection*{The reverse numerical range}
For $S\in\pB(\eX)$ and $x\in\SX$ let $W_x(S)=\{\langle Sx,\xi\rangle;\;\xi\in J(x)\}$. Hence
$W(S)=\bigcup_{x\in\SX}W_x(S)$. A number $\lambda$ belongs to $W_x(S)$ if and only if there exists $\xi\in J(x)$ such that
$\langle(S-\lambda)x,\xi\rangle=0$, that is, by Lemma \ref{lem22}, if and only if $x\bj(S-\lambda)x$. If we reverse the orientation 
in this last relation, then we obtain the following notion. 

\begin{definition} \label{def26}
The {\em local reverse numerical range} of $S\in\pB(\eX)$ at $x\in\SX$ is
$$\Wr_x(S)=\big\{\lambda\in\bC;\;(S-\lambda)x\bj x\big\}
=\big\{\lambda\in\bC;\;\|Sx-\lambda x\|\leq\|Sx-\mu x\|\;\text{for all}\;\mu\in\bC\big\}.$$
The {\em reverse numerical range} of $S$ is $\Wr(S)=\bigcup_{x\in\SX}\Wr_x(S)$.
\end{definition}

\noindent 
Thus the two ranges are
\begin{equation} \label{eq22}
W_x(S)=\big\{\lambda\in\bC;\;x\bj(S-\lambda)x\big\},\qquad \Wr_x(S)=\big\{\lambda\in\bC;\;(S-\lambda)x\bj x\big\}.
\end{equation}
The numerical range and the reverse numerical range are defined by the two orientations of Birkhoff--James
orthogonality, applied to the same pair of vectors. For operators on $\bC^n$ with an arbitrary norm, the notion was
introduced by Knockaert \cite{Kno1} in connection with linear prediction.

A number $\lambda$ belongs to $\Wr_x(S)$ if and only if $\lambda x$ is a best approximation of $Sx$ in the
one-dimensional subspace $\bC x$. Hence $\Wr_x(S)$ is the set of minimisers of the convex function
$\mu\mapsto\|Sx-\mu x\|$. This function is continuous and tends to infinity as $|\mu|\to\infty$. Therefore $\Wr_x(S)$
is a non-empty compact convex set. If $\lambda\in\Wr_x(S)$, then $\|Sx-\lambda x\|\leq\|Sx\|$, and therefore
$|\lambda|\leq2\|Sx\|$. Hence $\Wr(S)\subseteq\overline{D(0,2\|S\|)}$ (see also \cite[Theorem 5]{Kno1}). In Section
\ref{sec04} we determine the best constant in this inclusion. Clearly $\Wr(\alpha S+\beta I)=\alpha\Wr(S)+\beta$ for
all $\alpha,\beta\in\bC$, and the eigenvalues are included in it, that is,  $\sigma_p(S)\subseteq\Wr(S)$. If $\eX$ is a Hilbert space, 
Birkhoff--James orthogonality is symmetric and therefore $\Wr_x(S)=W_x(S)=\{\langle Sx,x\rangle\}$ for every $x\in\SX$, which
gives $\Wr(S)=W(S)$.

The following lemma connects the reverse numerical range of $S$ with the numerical range of its resolvent. For
$\lambda=0$ it is a reformulation of \cite[Proposition 1.6]{Poo}.

\begin{lemma} \label{lem27}
Let $S\in\pB(\eX)$. A number $\lambda\in\bC\setminus\sigma(S)$ belongs to $\Wr(S)$ if and only if 
$0\in W\big((S-\lambda)^{-1}\big)$. In particular, if $S$ is invertible, then $0\in\Wr(S)$ if and only if $0\in W(S^{-1})$.
\end{lemma}

\begin{proof}
Let $\lambda\in\bC\setminus\sigma(S)$. Suppose first that $\lambda\in\Wr(S)$. Then there exists $x\in\SX$ such that
$(S-\lambda)x\bj x$. Since $\lambda\notin\sigma(S)$, we have $(S-\lambda)x\ne 0$. Let
$y=\frac{1}{\|(S-\lambda)x\|}(S-\lambda)x\in\SX$. By homogeneity, $y\bj x$. By Lemma \ref{lem22}, there exists
$\eta\in J(y)$ such that $\langle x,\eta\rangle=0$. Since $(y,\eta)\in\PiX$ and
$\langle (S-\lambda)^{-1}y,\eta\rangle=\frac{1}{\|(S-\lambda)x\|}\langle x,\eta\rangle=0$, we have
$0\in W\big((S-\lambda)^{-1}\big)$.

Suppose now that $0\in W\big((S-\lambda)^{-1}\big)$. Then there exists $(x,\xi)\in\PiX$ such that
$\langle(S-\lambda)^{-1}x,\xi\rangle=0$. The vector $(S-\lambda)^{-1}x$ is non-zero because $x\neq0$. Let
$z=\frac{1}{\|(S-\lambda)^{-1}x\|}(S-\lambda)^{-1}x\in\SX$. Then
$\langle z,\xi\rangle=\frac{1}{\|(S-\lambda)^{-1}x\|}\langle (S-\lambda)^{-1}x,\xi\rangle=0$. It follows, by Lemma \ref{lem22}, that
$x\bj z$. Since $(S-\lambda)z=\frac{1}{\|(S-\lambda)^{-1}x\|}\,x$, homogeneity gives $(S-\lambda)z\bj z$, that is,
$\lambda\in\Wr(S)$.
\end{proof}

In the Corollary of \cite[Theorem 2]{Kno1}, the algebra numerical range appears in place of the spatial numerical range in Lemma
\ref{lem27}. Since the spatial numerical range of an operator is contained in its algebra numerical range, Lemma \ref{lem27} gives the 
implication $\lambda\in\Wr(S)\Rightarrow 0\in V\big((S-\lambda)^{-1}\big)$. The reverse implication requires that 
$0\in V\big((S-\lambda)^{-1}\big)$ implies $0\in W\big((S-\lambda)^{-1}\big)$. This holds if $W\big((S-\lambda)^{-1}\big)$ is closed 
and convex, for instance in a finite-dimensional Hilbert space (see \cite{Kno2}). However, as the following example shows, in general 
it fails. 

\begin{example} \label{ex28}
Let $\eX=\ell^1(2)$ and let 
$$Q=\begin{pmatrix}2&-i\\1&-2+2i\end{pmatrix}\in\pB(\eX).$$
The matrix $Q$ is invertible since $\det(Q)=-4+5i\ne 0$.
Let $x=(x_1,x_2)\in\SX$ and $\xi=(\xi_1,\xi_2)\in J(x)$. By \cite[Lemma 3.1]{BD}, $\xi_j=\frac{\bar x_j}{|x_j|}$ if $x_j\neq0$, and $|\xi_j|\leq1$ 
if $x_j=0$. If $x_2=0$, then $|x_1|=1$ and $\langle Qx,\xi\rangle=2+x_1\xi_2$. Similarly, if $x_1=0$, then $|x_2|=1$ and 
$\langle Qx,\xi\rangle=-2+2i-ix_2\xi_1$. In both cases $\langle Qx,\xi\rangle\neq0$ since $|x_1\xi_2|\leq1$ and $|x_2\xi_1|\leq1$. 
Suppose that $x_1\neq0$ and $x_2\neq0$. Then
$$\langle Qx,\xi\rangle=|x_1|(2+w)+|x_2|(-2+2i-i\bar w),\qquad\text{where}\quad w=\frac{x_1\bar x_2}{|x_1||x_2|}.$$
If this number were $0$, then $(2+w)\overline{(-2+2i-i\bar w)}$ would be a negative real number. However, this product
is equal to $-4-4i-2w+iw^2$, and its imaginary part $-4-2\im w+\re(w^2)$ is at most $-1$.
Hence $0\notin W(Q)$. On the other hand, for $x=(1,0)$ and $\xi=\big(1,-\frac{2+4i}{5}\big)$, and for $x=(0,1)$ and
$\xi=\big(\frac{4-2i}{5},1\big)$, we obtain the numbers $\mu_1=\frac{8-4i}{5}$ and $\mu_2=\frac{-12+6i}{5}$ in $W(Q)$.
Since $\frac35\mu_1+\frac25\mu_2=0$, we have $0\in V(Q)$. 

Thus, for the operator $S=Q^{-1}$, we have $0\in V(S^{-1})$, but $0\notin\Wr(S)$ by Lemma \ref{lem27}. Hence the Corollary of 
\cite[Theorem 2]{Kno1} does not hold for arbitrary norms.\qed
\end{example}

\subsection*{Dimension two}
The theorem of James and Day applies only to spaces of dimension at least $3$. There exist two-dimensional complex normed spaces 
with symmetric Birkhoff--James orthogonality which are not inner product spaces (see \cite[Chapter III]{Oma}). Therefore we need the 
following lemma in the proof of Theorem \ref{theo34}.

\begin{lemma} \label{lem29}
Let $\eX$ be a two-dimensional complex normed space. If $\eX$ has the \soi property, then $\pB(\eX)$ has the algebraic \soi property.
\end{lemma}

\begin{proof}
Assume that $\eX$ has the \soi property. First, we will show that $V(S)=-V(S)$ for every traceless operator $S\in\pB(\eX)$.
Thus, let $S$ be traceless. If $S=0$, there is nothing to prove. Suppose that $S\ne 0$. By the Cayley--Hamilton theorem,
$S^2=\mu^2I$ for some $\mu\in\bC$, and $\sigma(S)=\{\mu,-\mu\}$. Let $\kappa\in\bC\setminus\{\mu,-\mu\}$ and let
$T=\kappa I+S$. Then $T$ is invertible, and $T^2-2\kappa T+(\kappa^2-\mu^2)I=0$, which gives
$$T^{-1}=\frac{1}{\kappa^2-\mu^2}\,(\kappa I-S).$$
We see that $0\in W(T)$ if and only if $-\kappa\in W(S)$, and, similarly, $0\in W(T^{-1})$ if and only if $\kappa\in W(S)$. 
By the \soi property, it follows that $-\kappa\in W(S)$ if and only if $\kappa\in W(S)$, for all
$\kappa\ne\pm\mu$. Since $\pm\mu\in\sigma(S)\subseteq W(S)$, we conclude that $W(S)=-W(S)$. In finite dimensions $W(S)$ is
compact, and therefore $V(S)=\conv W(S)$ by \cite[Theorem 6]{SW}. Hence $V(S)=\conv W(S)=\conv\big(- W(S)\big)=-V(S)$.

Let now $T\in\pB(\eX)$ be an arbitrary invertible operator. If $T$ is a scalar multiple of $I$, then $0\notin V(T)$ and $0\notin V(T^{-1})$.
Otherwise, let $\kappa=\frac12\tr(T)$ and let $S=T-\kappa I$. Then $S\ne0$ and $\tr(S)=0$. As above, $S^2=\mu^2I$ with
$\kappa\neq\pm\mu$, since $\kappa\pm\mu$ are the eigenvalues of $T$. Moreover, $T^{-1}=\frac{1}{\kappa^2-\mu^2}(\kappa I-S)$. 
If $0\in V(T)$, then $-\kappa\in V(S)$. Since $V(S)=-V(S)$, we have $\kappa\in V(S)$, which implies $0\in V(T^{-1})$. Interchanging $T$
and $T^{-1}$ gives that $0\in V(T^{-1})$ implies $0\in V(T)$. Thus, $\pB(\eX)$ has the algebraic \soi property.
\end{proof}

\section{Characterisations of the symmetry of Birkhoff--James orthogonality} \label{sec03}

In this section we prove the main result of the paper. We begin with auxiliary results. As usual, for non-zero $u\in\eX$ and 
$\eta\in\eX^*$, let $u\otimes\eta$ be the rank-one operator $x\mapsto\langle x,\eta\rangle u$. Rank-one perturbations of the
identity were used for the same purpose by Ghosh \cite{Gho}.

\begin{proposition} \label{prop31}
Let $(u,\eta)\in\PiX$ and $\nu\in\bC$ with $0<|\nu|<1$. For the operator $T=I-\nu\,u\otimes\eta$, the following assertions hold.
\begin{enumerate}
\item $\re\big(\langle Tx,\xi\rangle\big)\geq1-|\nu|$ for all $(x,\xi)\in\PiX$. In particular, $T$ satisfies
the half-plane condition, and $0\notin V(T)$.
\item $T$ is invertible, and its inverse is $T^{-1}=I+\frac{\nu}{1-\nu}\,u\otimes\eta$.
\item If $(x,\xi)\in\PiX$, then $\langle T^{-1}x,\xi\rangle=0$ if and only if $\langle x,\eta\rangle\langle u,\xi\rangle=1-\frac{1}{\nu}$.
\end{enumerate}
\end{proposition}

\begin{proof}
(i) If $(x,\xi)\in\PiX$, then $\langle Tx,\xi\rangle=1-\nu\langle x,\eta\rangle\langle u,\xi\rangle$. Since
$|\langle x,\eta\rangle|\leq1$ and $|\langle u,\xi\rangle|\leq1$, we get $\re\big(\langle Tx,\xi\rangle\big)\geq1-|\nu|$.

(ii) Since $\|u\otimes\eta\|=\|u\|\|\eta\|=1$, we have $\|I-T\|=|\nu|<1$. Hence the operator $T$ is invertible. Because of 
$\langle u,\eta\rangle=1$ one has $(u\otimes\eta)^2=u\otimes\eta$. Now a straightforward computation shows that 
$T^{-1}=I+\frac{\nu}{1-\nu}\,u\otimes\eta$.

(iii) Since $T^{-1}=I+\frac{\nu}{1-\nu}\,u\otimes\eta$, by (ii), we have 
$\langle T^{-1}x,\xi\rangle=1+\frac{\nu}{1-\nu}\langle x,\eta\rangle\langle u,\xi\rangle$ for every $(x,\xi)\in\PiX$. 
It is clear that $1+\frac{\nu}{1-\nu}\langle x,\eta\rangle\langle u,\xi\rangle=0$ if and only if 
$\langle x,\eta\rangle\langle u,\xi\rangle=1-\frac{1}{\nu}$.
\end{proof}

Let
$$\Lambda(\eX)=\big\{\langle x,\eta\rangle\langle u,\xi\rangle;\;(x,\xi),(u,\eta)\in\PiX\big\}.$$
Note that $1\in\Lambda(\eX)$ and $\Lambda(\eX)\subseteq\overline{\bD}$, the closed unit disc.
In the case of a Hilbert space $\eH$, we have $J(x)=\{x\}$ for every $x\in\SH$. If $x,u\in\SH$, then 
$\langle x,u\rangle\langle u,x\rangle=|\langle x,u\rangle|^2\in[0,1]$, which means that $\Lambda(\eH)=[0,1]\subseteq\overline{D(1,1)}$.

\begin{corollary} \label{cor32}
If there exists a number $\omega\in\Lambda(\eX)$ such that $\omega\textcolor{blue}{\notin}\overline{D(1,1)}$, then $\eX$ does not have the 
\soi property and $\pB(\eX)$ does not have the algebraic \soi property. More precisely, if $(x,\xi),(u,\eta)\in\PiX$ are such that
\begin{equation} \label{eq31}
\big|1-\langle x,\eta\rangle\langle u,\xi\rangle\big|>1,
\end{equation}
and if $\nu=\big(1-\langle x,\eta\rangle\langle u,\xi\rangle\big)^{-1}$ and $T=I-\nu\,u\otimes\eta$, then $T$ satisfies the half-plane 
condition, and $0\in W(T^{-1})$. Consequently, the \soi property fails locally at $T$, and the algebraic \soi property of $\pB(\eX)$
fails at $T$.
\end{corollary}

\begin{proof}
We have $\langle x,\eta\rangle\langle u,\xi\rangle=1-\nu^{-1}$ and $0<|\nu|<1$ by \eqref{eq31}.
Thus, by Proposition \ref{prop31},  $T$ satisfies the half-plane condition, and $0\in W(T^{-1})$.
Note that $0\in W(T^{-1})\subseteq V(T^{-1})$, while $0\notin V(T)\supseteq\overline{W(T)}$.
\end{proof}

The next lemma is the key step in the proof of our main theorem.

\begin{lemma} \label{lem33}
If Birkhoff--James orthogonality in $\eX$ is not symmetric, then there exist $(x,\xi),(u,\eta)\in\PiX$ such that 
$\big|1-\langle x,\eta\rangle\langle u,\xi\rangle\big|>1$.
\end{lemma}

\begin{proof}
Since Birkhoff--James orthogonality is not symmetric, there exist vectors $u,y\in\SX$ such that $u\bj y$ but $y\not\bj u$. 
By Lemma \ref{lem22}, there exists $\eta\in J(u)$ such that $\langle y,\eta\rangle=0$, and the same lemma implies
$\langle u,\xi\rangle\ne0$ for every $\xi\in J(y)$. Let $\Omega=\{\langle u,\xi\rangle;\;\xi\in J(y)\}$.
Since the set $J(y)$ is weak$^*$ compact and convex, and the mapping $\xi\mapsto\langle u,\xi\rangle$ is weak$^*$ continuous and
linear, it follows that $\Omega$ is a compact convex subset of $\bC$ which does not contain $0$. Hence there exist $\theta\in[0,2\pi)$ 
and $\delta>0$ such that
\begin{equation} \label{eq32}
\re\big(e^{i\theta}\lambda\big)\leq-\delta\qquad\text{for all}\;\lambda\in\Omega .
\end{equation}

For $0<s<1$, let $y_s=y+se^{i\theta}u$. Each $y_s$ is a non-zero vector since $\|y_s\|\geq 1-s>0$. Let $x_s=\frac{1}{\|y_s\|}y_s$, 
and choose $\xi_s\in J(x_s)$. Thus $(x_s,\xi_s)\in\PiX$. We claim that there is $s_0\in(0,1)$ such that
\begin{equation} \label{eq33}
\re\big(e^{i\theta}\langle u,\xi_s\rangle\big)\leq-\tfrac12\delta\qquad\text{for all}\;s\in(0,s_0).
\end{equation}
Suppose that this is not true. Then there exists a sequence $(s_n)_{n=1}^{\infty}\subseteq (0,1)$ such that $s_n\to0$ and
$\re\big(e^{i\theta}\langle u,\xi_{s_n}\rangle\big)>-\frac12\delta$ for all $n$. Since the closed unit ball of $\eX^*$ is
weak$^*$ compact, the sequence $(\xi_{s_n})_{n=1}^\infty$ has a weak$^*$ cluster point $\xi_0$, and $\|\xi_0\|\leq1$. 
Since $y_s=\|y_s\|x_s$, we have $\|x_s-y_s\|=\big|1-\|y_s\|\big|\leq\|y-y_s\|=s$. Hence
$\|x_s-y\|\leq\|x_s-y_s\|+\|y_s-y\|\leq2s$, and therefore $\|x_{s_n}-y\|\to0$. Now we have
$|\langle y,\xi_{s_n}\rangle-1|=|\langle y-x_{s_n},\xi_{s_n}\rangle|\leq\|y-x_{s_n}\|\to0$
which implies $\langle y,\xi_0\rangle=1$. Thus $\xi_0\in J(y)$ and therefore $\langle u,\xi_0\rangle\in\Omega$. 
On the other hand, $\langle u,\xi_0\rangle$ is a cluster point of the sequence $(\langle u,\xi_{s_n}\rangle)_{n=1}^\infty$, and therefore 
$\re\big(e^{i\theta}\langle u,\xi_0\rangle\big)\geq-\frac12\delta$. This contradicts \eqref{eq32}. Hence there exists $s_0\in(0,1)$ such 
that \eqref{eq33} holds.

Let $s\in(0,s_0)$ and denote $\lambda_s=\langle x_s,\eta\rangle\langle u,\xi_s\rangle$. Since $\langle y,\eta\rangle=0$ and 
$\langle u,\eta\rangle=1$, we get
$$\lambda_s=\frac{1}{\|y_s\|}\langle y+se^{i\theta}u,\eta\rangle\langle u,\xi_s\rangle=\frac{s}{\|y_s\|}e^{i\theta}\langle u,\xi_s\rangle.$$
By \eqref{eq33}, $\re(\lambda_s)=\frac{s}{\|y_s\|}\re\big(e^{i\theta}\langle u,\xi_s\rangle\big)\leq-\frac{s\delta}{2\|y_s\|}<0$. Hence
$|1-\lambda_s|^2=1-2\re(\lambda_s)+|\lambda_s|^2>1$.
Thus we may choose the pairs $(x_s,\xi_s), (u,\eta)\in\PiX$ to fulfil the conclusion of this lemma.
\end{proof}

We are now ready to prove the main result. It shows that the symmetry of Birkhoff--James orthogonality, which is a property of the 
norm, can be recognised by numerical ranges of operators on $\eX$. Moreover, it gives a simple condition in terms of supporting 
functionals.

\begin{theorem} \label{theo34}
For a complex Banach space $\eX$ of dimension at least $2$, the following assertions are equivalent.
\begin{enumerate}
\renewcommand{\theenumi}{\alph{enumi}}\renewcommand{\labelenumi}{(\theenumi)}
\item Birkhoff--James orthogonality in $\eX$ is symmetric.
\item $W(S)=\Wr(S)$ for every $S\in\pB(\eX)$.
\item $\eX$ has the \soi property.
\item $\pB(\eX)$ has the algebraic \soi property.
\item $\Lambda(\eX)\subseteq\overline{D(1,1)}$.
\end{enumerate}
\end{theorem}

\begin{proof}
(a)$\Rightarrow$(b). If Birkhoff--James orthogonality is symmetric, then, by \eqref{eq22},
$W_x(S)=\Wr_x(S)$ for every $x\in\SX$. Taking the union over $x\in\SX$ gives (b).

(a)$\Rightarrow$(c) and (a)$\Rightarrow$(d). Suppose first that $\dim\eX\geq3$. By the theorem of James and Day
\cite{Jam47a,Day}, which holds also for complex spaces \cite{Boh,AGKRZ}, $\eX$ is a Hilbert space. Hence (c) and (d)
hold by Lemma \ref{lem25}. Suppose now that $\dim\eX=2$. Then $W(T)$ is compact for every $T\in\pB(\eX)$. Hence (c)
follows from Lemma \ref{lem23}, and (d) follows from (c) by Lemma \ref{lem29}.

(b)$\Rightarrow$(e), (c)$\Rightarrow$(e) and (d)$\Rightarrow$(e). Suppose that (e) does not hold. Then there exist
$(x,\xi),(u,\eta)\in\PiX$ which satisfy \eqref{eq31}. By Corollary \ref{cor32}, there exists $T\in\pB(\eX)$ which satisfies
the half-plane condition and $0\in W(T^{-1})$. In particular, neither (c) nor (d) holds. Moreover, $0\notin W(T)$, while
$0\in\Wr(T)$ by Lemma \ref{lem27}. Hence (b) does not hold.

(e)$\Rightarrow$(a). This is the contrapositive of Lemma \ref{lem33}.
\end{proof}

The implication (c)$\Rightarrow$(a) answers a question of Poon \cite{Poo} and confirms a conjecture of Ghosh \cite{Gho}.
Ghosh proved it for smooth and almost locally uniformly rotund spaces \cite[Theorem 1.5]{Gho}, for finite-dimensional smooth
or strictly convex spaces \cite[Corollary 3.4]{Gho}, for real polyhedral spaces \cite[Theorem 3.7]{Gho}, and for real planes
\cite[Theorem 3.10]{Gho}.

\begin{corollary} \label{cor35}
For a complex Banach space $\eX$ of dimension at least $3$, the following assertions are equivalent.
\begin{enumerate}
\renewcommand{\theenumi}{\alph{enumi}}\renewcommand{\labelenumi}{(\theenumi)}
\item $\eX$ is a Hilbert space.
\item $\eX$ has the \soi property.
\item $\pB(\eX)$ has the algebraic \soi property.
\end{enumerate}
\end{corollary}

\begin{proof}
This follows from Theorem \ref{theo34} and the theorem of James and Day.
\end{proof}

\noindent 
Under additional smoothness or rotundity assumptions, the equivalence of (a) and (b) was proved by Ghosh \cite[Corollary 3.4]{Gho}.
In dimension two, Theorem \ref{theo34} characterises the \soi property by the symmetry of Birkhoff--James orthogonality. The examples 
in \cite[Chapter III]{Oma} and \cite{Poo} show that this class of planes is larger than the class of two-dimensional Hilbert spaces.

\begin{corollary} \label{cor36}
Let $\eX$ be a complex Banach space of dimension at least $3$. Then $\eX$ has the \soi property if and only if $\eX^*$ has the \soi property.
\end{corollary}

\begin{proof}
If $\eX$ is a Hilbert space, then $\eX^*$ is a Hilbert space as well. Conversely, if $\eX^*$ is a Hilbert space, then $\eX^{**}$ is a
Hilbert space, and $\eX$ is isometrically embedded in $\eX^{**}$. Hence $\eX$ is a Hilbert space. Apply Corollary \ref{cor35}.
\end{proof}

\begin{remark} \label{rem37}
(i) The part (b)$\Rightarrow$(e) of the proof of Theorem \ref{theo34} shows that (b) may be replaced by the weaker assertion
$$0\in W(T)\iff0\in\Wr(T)\qquad\text{for every invertible}\;T\in\pB(\eX).$$
By Lemma \ref{lem27}, this is the \soi property without closures. Assertion (b) may also be replaced by
$W(S)\setminus\sigma(S)=\Wr(S)\setminus\sigma(S)$ for every $S\in\pB(\eX)$. In contrast to assertion (c), no closures are involved.

(ii) If $\eX$ is a complex Banach space without the \soi property, then there are $(u,\eta)\in\PiX$ and $\nu\in\bC$ with 
$0<|\nu|<1$ such that the operator $T=I-\nu\,u\otimes\eta$ satisfies the half-plane condition and $0\in W(T^{-1})$. We have 
shown that these operators are crucial: a complex Banach space $\eX$ of dimension at least $2$ has the \soi property if and only if the 
\soi property holds locally at every operator $I-\nu\,u\otimes\eta$ with $(u,\eta)\in\PiX$ and $0<|\nu|<1$.

(iii) The proofs of Proposition \ref{prop31} and Lemma \ref{lem33} work for real Banach spaces, with $\theta\in\{0,\pi\}$ and real $\nu$. 
Together with the real versions of Lemma \ref{lem23} and Lemma \ref{lem25} and the theorem of James and Day, they show that a
real Banach space has the \soi property if and only if Birkhoff--James orthogonality is symmetric in it. For real
planes this is \cite[Theorem 3.10]{Gho}. In higher dimensions, Ghosh needed smoothness or rotundity assumptions
\cite[Corollary 3.4]{Gho}, while our argument needs none. In the real case we have $\Lambda(\eX)\subseteq[-1,1]$, and assertion (e) 
means $\Lambda(\eX)\subseteq[0,1]$. In this form, the equivalence (a)$\Leftrightarrow$(e) is known for real spaces. Indeed, 
for $x,u\in\SX$, Lemma \ref{lem22} shows that $x\bj (ax+u)$ if and only if $a=-\langle u,\xi\rangle$ for some $\xi\in J(x)$. James 
\cite[Theorem 2.4]{Jam47b} proved that $ab\geq0$ whenever $x\bj (ax+u)$ and $u\bj (bu+x)$, provided Birkhoff--James orthogonality 
is symmetric. Martini and Swanepoel \cite[Corollary 3]{MS} proved the converse for real planes. It extends to all real spaces, since 
$x\bj y$ depends only on the plane spanned by $x$ and $y$ (see also \cite[Theorem 4.10]{AMW}). Our Theorem \ref{theo34} extends 
this equivalence to complex spaces, where the condition takes the form $|1-\lambda|\leq1$.
\end{remark}

\section{The asymmetry constant of Birkhoff--James orthogonality} \label{sec04}

Theorem \ref{theo34} is qualitative; it says when Birkhoff--James orthogonality is symmetric. However, it does not say how far 
Birkhoff--James orthogonality is from being symmetric when it is not. In this section we measure the asymmetry by a single number. 
For $(x,\xi)\in\PiX$, let $P_{x,\xi}=I-x\otimes\xi$. It is easy to see that this is the non-zero projection onto the kernel of $\xi$ along $x$. 
Hence $\|P_{x,\xi}\|\geq 1$. By Lemma \ref{lem22}, $x\bj z$ for every $z\in\ker\xi$, and every pair $x\bj z$ with $x\in\SX$ arises in this 
way. The reverse relation $z\bj x$ means that $\|z\|\leq\|z+\mu x\|$ for all $\mu\in\bC$. Since $P_{x,\xi}(z+\mu x)=z$, this holds for
all $z\in\ker\xi$ if and only if $\|P_{x,\xi}\|=1$. Hence Birkhoff--James orthogonality is symmetric if and only if
$\|P_{x,\xi}\|=1$ for every $(x,\xi)\in\PiX$. This leads to the following definition. 

\begin{definition} \label{def41}
The {\em asymmetry constant} of Birkhoff--James orthogonality in $\eX$ is  
$$R(\eX)=\sup\big\{\|P_{x,\xi}\|;\;(x,\xi)\in\PiX\big\}.$$
\end{definition}

Since $x\otimes\xi$ is an idempotent of norm $1$ for every $(x,\xi)\in\PiX$, we see that $\|P_{x,\xi}\|\leq 2$. Hence, for every
complex Banach space of dimension at least $2$, the asymmetry constant is a number between $1$ and $2$. For later use, we
formulate these basic observations as a proposition.

\begin{proposition} \label{prop42}
For every complex Banach space $\eX$ of dimension at least $2$, the asymmetry constant $R(\eX)$ lies in $[1,2]$. It is equal to $1$ if 
and only if Birkhoff--James orthogonality in $\eX$ is symmetric.
\end{proposition}

\begin{proof}
This follows from the discussion before the proposition.
\end{proof}

For closed subspaces $\eM$ and $\eN$ such that $\eM\cap\eN=\{0\}$, let
$$\gamma(\eM,\eN)=\inf\{\dist(z,\eN);\;z\in\eM,\ \|z\|=1\}$$ 
be the minimum gap (see \cite[Chapter IV, \S4]{Kat}). In a Hilbert space, it is the sine of the minimal angle between $\eM$ and
$\eN$. If $P$ is a bounded projection with range $\eM$ and kernel $\eN$, then $\|P\|^{-1}=\gamma(\eM,\eN)$, 
by \cite[Chapter IV, (4.7)]{Kat}. Hence $R(\eX)^{-1}=\inf\big\{\gamma(\ker\xi,\bC x);\;(x,\xi)\in\PiX\big\}$.
The next lemma expresses this number directly in terms of Birkhoff--James orthogonality.

\begin{lemma} \label{lem43}
For the asymmetry constant of Birkhoff--James orthogonality in $\eX$ we have
$$R(\eX)^{-1}=\inf\big\{\|x+v\|;\;x\in\SX,\;v\in\eX,\;v\bj x\big\}.$$
\end{lemma}

\begin{proof}
Let $m$ denote the infimum. Let $x\in\SX$ and $v\in\eX$ be such that $v\bj x$. If $v=0$, then $\|x+v\|=1\geq R(\eX)^{-1}$. 
Suppose $v\ne0$ and denote $y=\frac{1}{\|v\|}v$. Then $y\bj x$, and by Lemma \ref{lem22} there exists $\eta\in J(y)$ such that
$\langle x,\eta\rangle=0$. Let $P=I-y\otimes\eta$. Since $\langle x+v,\eta\rangle=\|v\|$, we get
$P(x+v)=x+v-\|v\| y=x$. Hence $1=\|x\|\leq \|P\|\|x+v\|\leq R(\eX)\|x+v\|$ and therefore $\|x+v\|\geq R(\eX)^{-1}$. Since $x\in\SX$
and $v\bj x$ are arbitrary we conclude that $m\geq R(\eX)^{-1}$.

To prove the reverse inequality, let $(y,\eta)\in\PiX$ be arbitrary and let $P=I-y\otimes\eta$. If $z\in\eX$ is such that $Pz\ne0$, then 
$\langle Pz,\eta\rangle=\langle z-\langle z,\eta\rangle y,\eta\rangle=0$. It follows that $y\bj Pz$ by Lemma \ref{lem22}. 
Let $x=\frac{1}{\|Pz\|}Pz$ and $v=\frac{\langle z,\eta\rangle}{\|Pz\|}y$. Then $x\in\SX$ and $v\bj x$, since $y\bj Pz$ and Birkhoff--James 
orthogonality is homogeneous in both variables. Note that $x+v=\frac{1}{\|Pz\|}\big(Pz+\langle z,\eta\rangle y\big)=\frac{1}{\|Pz\|}z$.
Hence $\frac{\|Pz\|}{\|z\|}=\frac{1}{\|x+v\|}\leq\frac{1}{m}$ which gives $\|P\|\leq \frac{1}{m}$ for every pair $(y,\eta)\in\PiX$. It follows that 
$R(\eX)\leq \frac{1}{m}$, that is, $m\leq R(\eX)^{-1}$.
\end{proof}

For real normed spaces, Ji and Wu \cite{JW} introduced the constant
$$D(\eX)=\inf\Big\{\inf_{\lambda\in\bR}\|x+\lambda y\|;\;x,y\in\SX,\;\|x+y\|=\|x-y\|\Big\},$$
which measures the difference between Birkhoff--James orthogonality and isosceles orthogonality. They showed that $D(\eX)=1$ if and 
only if $\eX$ is an inner product space, and computed $D(\ell^p(2))$ \cite[Theorems 1 and 3]{JW}. By Lemma \ref{lem43}, the constant 
$R(\eX)^{-1}$ is obtained in the same way, with isosceles orthogonality replaced by the reverse Birkhoff--James orthogonality, that is,
$R(\eX)^{-1}=\inf\Big\{\inf_{\lambda\in\bC}\|x+\lambda y\|;\;x,y\in\SX,\;y\bj x\Big\}$. Other constants measuring differences between
orthogonality types in real normed spaces were studied, for instance, in \cite{PW}.

We now show that, in every Banach space, the \soi property holds locally at all operators which are sufficiently
close to the identity. The estimate is elementary, and it actually holds in every unital Banach algebra.

\begin{lemma} \label{lem44}
Let $\pA$ be a unital Banach algebra. If $a\in\pA$ is such that $\|e-a\|=r<1$, then $a$ is invertible and $\re(\lambda)\geq1-r$ for every
$\lambda\in V(a)$. If $r<\frac12$, then $\re(\lambda)\geq\frac{1-2r}{1-r}$ for every $\lambda\in V(a^{-1})$. In particular, 
if $T\in\pB(\eX)$ and $\|I-T\|<\frac12$, then $0\notin\overline{W(T)}$ and
$0\notin\overline{W(T^{-1})}$. Hence the \soi property holds locally at $T$.
\end{lemma}

\begin{proof}
It is well known that $a$ is invertible if $\|e-a\|=r<1$. Let $\varphi$ be a state on $\pA$. Then $\varphi(a)=1-\varphi(e-a)$ and 
$|\varphi(e-a)|\leq r$. Hence $\re\big(\varphi(a)\big)=1-\re\big(\varphi(e-a)\big)\geq 1-|\varphi(e-a)|\geq1-r$. Since 
$e-a^{-1}=a^{-1}(a-e)$, we have $\|e-a^{-1}\|\leq r\|a^{-1}\|$. Hence $1=\|e\|\geq\|a^{-1}\|-\|e-a^{-1}\|\geq\|a^{-1}\|(1-r)$, that is, 
$\|a^{-1}\|\leq\frac{1}{1-r}$, and therefore $\|e-a^{-1}\|\leq\frac{r}{1-r}$. If $r<\frac12$, then $\frac{r}{1-r}<1$, and
the first part applied to $a^{-1}$ gives $\re(\lambda)\geq1-\frac{r}{1-r}=\frac{1-2r}{1-r}$ for every $\lambda\in V(a^{-1})$.
The last assertions follow from $\overline{W(T)}\subseteq V(T)$.
\end{proof}

By Lemma \ref{lem44}, in every Banach space the \soi property holds locally at every operator in the open ball of radius $\frac12$ 
around $I$. This leads to the following definition. 

\begin{definition} \label{def45}
The {\em \soi radius} of $\eX$ is
$$s_W(\eX)=\sup\big\{r\in(0,1];\;\text{the \soi property holds locally at every}\;T\;\text{with}\;\|I-T\|<r\big\}.$$
\end{definition}

By the definition and Lemma \ref{lem44}, $\frac12\leq s_W(\eX)\leq1$. If $\eX$ has the \soi property, then $s_W(\eX)=1$.
The next theorem shows that the asymmetry constant determines the \soi radius. Together with Proposition \ref{prop42} and 
Theorem \ref{theo34}, this shows that the converse also holds.

\begin{theorem} \label{theo46}
For the asymmetry constant $R(\eX)$, the following hold.
\begin{enumerate}
\item $s_W(\eX)=R(\eX)^{-1}$. 
\item The constant $R(\eX)$ is the best possible number $c>0$ such that $\Wr(S)\subseteq\overline{D(0,c\|S\|)}$ for every
$S\in\pB(\eX)$, that is,
$$R(\eX)=\sup\big\{|\lambda|;\;\lambda\in\Wr(S),\;S\in\pB(\eX),\;\|S\|\leq1\big\}.$$
Moreover, $|\lambda|\leq R(\eX)\|Sx\|$ for every $x\in\SX$ and every $\lambda\in\Wr_x(S)$.
\end{enumerate}
\end{theorem}

\begin{proof}
(i) Let $T=I-C$ with $\|C\|=r<R(\eX)^{-1}$. By the first part of Lemma \ref{lem44}, we have $0\textcolor{blue}{\notin}\overline{W(T)}$. Let 
$\lambda\in W(T^{-1})$ be arbitrary and let $(x,\xi)\in\PiX$ be such that $\lambda=\langle T^{-1}x,\xi\rangle$.  Denote 
$P=I-x\otimes\xi$. Then $Px=x-\langle x,\xi\rangle x=0$ and $PT^{-1}x=T^{-1}x-\lambda x$. Since $T^{-1}x-CT^{-1}x=x$, we have 
$PT^{-1}x=PCT^{-1}x$. We conclude that $T^{-1}x=\lambda x+PCT^{-1}x$. It follows that
\begin{equation} \label{eq41}
\|T^{-1}x\|=\|\lambda x+PCT^{-1}x\|\leq |\lambda|+\|P\|\|C\|\|T^{-1}x\|\leq |\lambda|+R(\eX) r\| T^{-1}x\|.
\end{equation}
Note that $\|T\|=\|I-C\|\leq 1+r$, that is, $\|T\|^{-1}\geq \frac{1}{1+r}$, and $\|T\|^{-1}\leq\| T^{-1}x\|$. Hence 
$\|T^{-1}x\|\geq \frac{1}{1+r}$.  Using this in \eqref{eq41} gives
$$|\lambda|\geq\|T^{-1}x\|\big(1-R(\eX)r\big)\geq \frac{1-R(\eX)r}{1+r}>0.$$
Hence $W(T^{-1})$ is disjoint from the open disc $D\big(0,\frac{1-R(\eX)r}{1+r}\big)$. In particular,
$0\textcolor{blue}{\notin}\overline{W(T^{-1})}$. We conclude that the \soi property holds locally at $T$. Hence $s_W(\eX)\geq R(\eX)^{-1}$.

If $R(\eX)=1$, then $s_W(\eX)=1$, by the previous paragraph. Suppose therefore that $R(\eX)>1$ and let $r>0$ be such that
$r^{-1}<R(\eX)$. Hence there exists $(x,\xi)\in\PiX$ such that $\| P_{x,\xi}\|=\|I-x\otimes\xi\|>\max\{1, r^{-1}\}$. It follows that
there exist $u\in\SX$ and $\eta\in\SXd$ such that $|\langle P_{x,\xi}u,\eta\rangle|>\max\{1, r^{-1}\}$. Denote
$\omega=\langle P_{x,\xi}u,\eta\rangle$. Since $|\langle u,\eta\rangle|\leq 1$, we see that $\omega\ne\langle u,\eta\rangle$.
Let $T=I-\frac{1}{\omega}u\otimes \eta$. Then $\| I-T\|=\frac{1}{|\omega|}<\min\{1,r\}$, and therefore $T$ is invertible. It is easy to 
see that $T^{-1}=I+\frac{1}{\omega-\langle u,\eta\rangle}u\otimes\eta$. Since 
$\omega=\langle u,\eta\rangle-\langle u,\xi\rangle\langle x,\eta\rangle$, we get
$$\langle T^{-1}x,\xi\rangle=1+\frac{\langle x,\eta\rangle\langle u,\xi\rangle}{\omega-\langle u,\eta\rangle}=
\frac{\langle u,\eta\rangle-\langle u,\xi\rangle\langle x,\eta\rangle-\langle u,\eta\rangle+\langle x,\eta\rangle\langle u,\xi\rangle}
{\omega-\langle u,\eta\rangle}=0 .$$
Hence $0\in W(T^{-1})$. Since $0\notin\overline{W(T)}$ by Lemma \ref{lem44}, the \soi property fails locally at $T$. Therefore
$s_W(\eX)\leq R(\eX)^{-1}$.

(ii) Let $x\in\SX$ and $\lambda\in\Wr_x(S)$. If $\lambda=0$, there is nothing to prove. If $\lambda\ne0$, then $(S-\lambda)x\bj x$ 
and hence $\lambda^{-1}(S-\lambda)x\bj x$. By Lemma \ref{lem43}, $\|x+\lambda^{-1}(S-\lambda)x\|\geq R(\eX)^{-1}$. Since 
$Sx=\lambda\big(x+\lambda^{-1}(S-\lambda)x\big)$, we get $\|Sx\|\geq|\lambda|R(\eX)^{-1}$. This is the last assertion of (ii). Since 
$\Wr(S)=\bigcup_{x\in\SX}\Wr_x(S)$ and $\|Sx\|\leq\|S\|$, the supremum in (ii) is at most $R(\eX)$. 

Let $\varepsilon>0$. By Lemma \ref{lem43}, there exist $x\in\SX$ and $v\in\eX$ such that $v\bj x$ and 
$\|x+v\|^{-1}>R(\eX)-\varepsilon$. Note that $\|x+v\|\geq R(\eX)^{-1}>0$. Let $\xi\in J(x)$ and $S=\|x+v\|^{-1}(x+v)\otimes\xi$. 
Then $\|S\|=1$ and $(S-\|x+v\|^{-1})x=\|x+v\|^{-1}(x+v)-\|x+v\|^{-1}x=\|x+v\|^{-1}v$. It follows from $v\bj x$ that 
$\big(S-\|x+v\|^{-1}\big)x\bj x$, that is, $\|x+v\|^{-1}\in\Wr(S)$. Since $\|x+v\|^{-1}>R(\eX)-\varepsilon$, we conclude that the 
supremum in (ii) is at least $R(\eX)$.
\end{proof}

If $(x,\xi)\in\PiX$, $y\in\eX$, and $S=y\otimes\xi$, then $\Wr_x(S)$ is the set of all numbers $\lambda$ for which $\lambda x$ is
a best approximation of $y$ in $\bC x$. Hence, by Theorem \ref{theo46}~(ii), $R(\eX)$ is the supremum of quotients
$\frac{\|m\|}{\|y\|}$, where $y\ne0$ and $m$ is a best approximation of $y$ in a one-dimensional subspace of $\eX$.
In particular, Birkhoff--James orthogonality is symmetric if and only if such best approximations never exceed $\|y\|$ (cf.\ \cite{Ami}).

The reverse numerical range also controls the spectrum. For a non-empty compact set $K\subseteq\bC$, let $\widehat{K}$
denote its polynomially convex hull. Recall that the approximate point spectrum of $T\in\pB(\eX)$ is the subset $\sigma_{ap}(T)$ of
$\sigma(T)$ consisting of all numbers $\lambda$ for which there is a sequence $(x_n)_{n=1}^{\infty}\subseteq\SX$
such that $\lim\limits_{n\to\infty}\|Tx_n-\lambda x_n\|=0$.

\begin{corollary} \label{cor47}
Let $S\in \pB(\eX)$ be arbitrary. Then
\begin{equation} \label{eq42}
\|(S-\lambda)x\|\geq R(\eX)^{-1} \dist\big(\lambda,\Wr(S)\big)\qquad\text{for all}\;\lambda\in\bC\;\text{and}\;
x\in\SX .
\end{equation}
Consequently, $\sigma_{ap}(S)\subseteq\overline{\Wr(S)}$, and $\sigma(S)\subseteq\widehat{\overline{\Wr(S)}}$.
Moreover, 
$$\|(S-\lambda)^{-1}\|\leq \frac{R(\eX)}{\dist(\lambda,\Wr(S))}\qquad \text{for every}\; 
\lambda\in \bC\setminus\widehat{\overline{\Wr(S)}}.$$
\end{corollary}

\begin{proof}
Let $x\in\SX$. If $\mu\in\Wr_x(S)$, then $\mu-\lambda\in\Wr_x(S-\lambda)$, and the last assertion of Theorem \ref{theo46}~(ii) gives
$|\mu-\lambda|\leq R(\eX)\|(S-\lambda)x\|$. Hence \eqref{eq42} holds.

Suppose that $\lambda\in\sigma_{ap}(S)$. Then the infimum over $x\in\SX$ of the left-hand side of \eqref{eq42} is $0$.
It follows that $\lambda\in\overline{\Wr(S)}$. The set $\overline{\Wr(S)}$ is compact by Theorem \ref{theo46}~(ii). Let $U$ be 
the unbounded component of $\bC\setminus \overline{\Wr(S)}$. This set is connected, and it contains points of the resolvent set of $S$, 
since it is unbounded. Since $\partial\sigma(S)\subseteq\sigma_{ap}(S)\subseteq \overline{\Wr(S)}$, the set $U$ does not meet the boundary of $\sigma(S)$. Hence $U$ does not meet $\sigma(S)$, which gives 
$\sigma(S)\subseteq\bC\setminus U=\widehat{\overline{\Wr(S)}}$. 
For $\lambda\in U$ the operator $S-\lambda$ is invertible, and the estimate for $\|(S-\lambda)^{-1}\|$ follows from the first assertion.
\end{proof}

\begin{remark} \label{rem48}
(i) For a Hilbert space $\eH$, the asymmetry constant is $1$ and $\Wr(S)=W(S)$. Hence, Corollary \ref{cor47} reduces to the classical 
estimate $\|(S-\lambda)x\|\geq\dist(\lambda,W(S))$.

(ii) Theorem \ref{theo46} gives a second proof of the implications (c)$\Rightarrow$(a) and (d)$\Rightarrow$(a) of 
Theorem \ref{theo34}, which avoids Lemma \ref{lem33}. Indeed, each of the properties (c) and (d) gives $s_W(\eX)=1$. 
For (d) note that $\|I-T\|<1$ implies $0\notin V(T)$, hence $0\notin V(T^{-1})\supseteq\overline{W(T^{-1})}$. Since also
$0\notin\overline{W(T)}\subseteq V(T)$, the \soi property holds locally at every $T$ with $\|I-T\|<1$. Hence $R(\eX)=1$ 
by Theorem \ref{theo46}~(i), and Birkhoff--James orthogonality is symmetric by Proposition \ref{prop42}.
\end{remark}

\begin{corollary} \label{cor49}
The following hold.
\begin{enumerate}
\item If $\eY\subseteq\eX$ is a closed subspace of dimension at least $2$, then $R(\eY)\leq R(\eX)$.
\item $R(\eX)\leq R(\eX^*)$, and the equality holds if $\eX$ is reflexive.
\item If $-1\in\overline{\Lambda(\eX)}$, then $R(\eX)=2$.
\end{enumerate}
\end{corollary}

\begin{proof}
(i) Let $(x,\xi)\in\Pi_{\eY}$. By the Hahn--Banach theorem, $\xi$ has an extension $\tilde\xi\in\eX^*$ with
$\|\tilde\xi\|=1$. Then $(x,\tilde\xi)\in\PiX$, and $I_{\eY}-x\otimes\xi$ is the restriction of $I-x\otimes\tilde\xi$
to $\eY$. Hence $\|I_{\eY}-x\otimes\xi\|\leq\|I-x\otimes\tilde\xi\|\leq R(\eX)$, and therefore $R(\eY)\leq R(\eX)$.

(ii) Let $(x,\xi)\in\PiX$, and let $\hat x\in\eX^{**}$ be the canonical image of $x$. Then $(\xi,\hat x)\in\Pi_{\eX^*}$,
and $(I-x\otimes\xi)^*=I-\xi\otimes\hat x$. Since an operator and its adjoint have the same norm,
$R(\eX)\leq R(\eX^*)$. If $\eX$ is reflexive, apply this to $\eX^*$.

(iii) For $(x,\xi),(u,\eta)\in\PiX$ we have $\langle P_{x,\xi}u,\eta\rangle=1-\langle x,\eta\rangle\langle u,\xi\rangle$. Hence
$R(\eX)\geq\sup\big\{|1-\lambda|;\;\lambda\in\Lambda(\eX)\big\}$. Since $-1\in\overline{\Lambda(\eX)}$, we get $R(\eX)\geq2$, and 
$R(\eX)\leq2$ by Proposition \ref{prop42}.
\end{proof}

Every two-dimensional space is reflexive. Hence, by Corollary \ref{cor49}~(ii), Proposition \ref{prop42} and Theorem \ref{theo34}, 
Corollary \ref{cor36} holds also when $\dim\eX=2$.

We will close the section with a few examples. The first gives an explicit Banach space $\eX$ with $R(\eX)=2$. Then we compute the 
asymmetry constant of the two-dimensional spaces $\ell^p(2)$.

\begin{example} \label{ex410}
Let $\eX=\ell^1(n)$, $n\geq2$, and let $(e_1,\ldots, e_n)$ be the standard basis. Denote the standard basis of $\eX^*=\ell^\infty(n)$
by $(e_{1}^{*},\ldots,e_{n}^{*})$.
For $\alpha,\beta\in\overline{\bD}$, the functionals $\eta=e_{1}^{*}+\alpha e_{2}^{*}$ and $\xi=\beta e_{1}^{*}+e_{2}^{*}$ satisfy 
$\eta\in J(e_1)$ and $\xi\in J(e_2)$. Moreover, $\langle e_2,\eta\rangle\langle e_1,\xi\rangle=\alpha\beta$. Hence 
$\Lambda(\ell^1(n))=\overline{\bD}$, and $R(\ell^1(n))=2$ by Corollary \ref{cor49}~(iii). The part (ii) of the same corollary gives
$R(\ell^\infty(n))=2$. That $\ell^1(n)$ and $\ell^\infty(n)$ do not have the \soi property was shown in
\cite[Example 3.7 and Remark 3.8]{BD}.

Let $\eX=\eY\oplus_1\eZ$ or $\eX=\eY\oplus_\infty\eZ$ with $\eY\ne\{0\}\ne\eZ$. Arguments as above show that 
$\Lambda(\eX)=\overline{\bD}$. This covers $\ell^1$, $\ell^\infty$, $c_0$, $L^1(\mu)$,
$L^\infty(\mu)$ and $M(K)$, whenever the dimension is at least $2$. For all these spaces $R(\eX)=2$. By Corollary \ref{cor49}~(i), 
the same holds for every space which contains an isometric copy of $\ell^1(2)$ or of $\ell^\infty(2)$. For instance, let $K$ be a compact 
Hausdorff space with at least two points. By Urysohn's lemma there are $f,g\in C(K)$ with $\|f\|=\|g\|=1$ and disjoint supports. Then
$\|\alpha f+\beta g\|=\max\{|\alpha|,|\beta|\}$, and $C(K)$ contains an isometric copy of $\ell^\infty(2)$. Hence $R(C(K))=2$.\qed
\end{example}

That $\ell^p(n)$, $1<p<\infty$, $p\ne2$, does not have the \soi property was conjectured in
\cite[Remark 3.8]{BD}. It was proved by Poon \cite[Theorem 1.9]{Poo} for $n=2$ and by Ghosh \cite[Corollary 3.4 and
Remark 3.5]{Gho} for $n\geq3$. The next proposition gives the asymmetry constant of $\ell^p(2)$.

\begin{proposition} \label{prop411}
Let $1<p<\infty$ and $q=\frac{p}{p-1}$. Then
$$R(\ell^p(2))=\max_{0\leq \tau\leq1}\big(\tau^{p-1}+(1-\tau)^{p-1}\big)^{1/p}\big(\tau^{q-1}+(1-\tau)^{q-1}\big)^{1/q}.$$
In particular, $R(\ell^p(2))=R(\ell^q(2))$. Moreover, $R(\ell^p(n))\geq R(\ell^p(2))$ for every $n\geq2$.
\end{proposition}

\begin{proof}
Let $U$ be a surjective linear isometry of $\eX=\ell^p(2)$. If $(x,\xi)\in\PiX$, then $(Ux,\xi\circ U^{-1})\in\PiX$
and $P_{Ux,\xi\circ U^{-1}}=UP_{x,\xi}U^{-1}$. Diagonal unitary matrices and the transposition of coordinates are
isometries of $\ell^p(2)$. Hence it suffices to consider $x=(s,t)$ with $s,t\geq0$ and $s^p+t^p=1$. Since the norm
is smooth, $J(x)=\{\xi\}$ with $\xi=(s^{p-1},t^{p-1})$. Let $w=(t^{p-1},-s^{p-1})$ and $\zeta=(t,-s)$. Then
$\langle w,\xi\rangle=0$, $\langle x,\zeta\rangle=0$ and $\langle w,\zeta\rangle=s^p+t^p=1$. Hence the operators
$P_{x,\xi}$ and $w\otimes\zeta$ both vanish at $x$ and both fix $w$. Since $x$ and $w$ form a basis of $\eX$, we
have $P_{x,\xi}=w\otimes\zeta$, and therefore $\|P_{x,\xi}\|=\|w\|_p\|\zeta\|_q$. Let $\tau=s^p$. Then
$\|w\|_p^p=(1-\tau)^{p-1}+\tau^{p-1}$ and, since $\frac qp=q-1$, $\|\zeta\|_q^q=(1-\tau)^{q-1}+\tau^{q-1}$. This gives the
formula. The right-hand side is symmetric in $p$ and $q$. Finally, $\ell^p(2)$ is isometric to a subspace of $\ell^p(n)$, and the last 
assertion follows from Corollary \ref{cor49}~(i).
\end{proof}

\begin{corollary} \label{cor412}
The function $p\mapsto R(\ell^p(2))$ is continuous on $(1,\infty)$. Moreover, $R(\ell^2(2))=1$ and
$$\lim_{p\to1}R(\ell^p(2))=\lim_{p\to\infty}R(\ell^p(2))=2.$$
Consequently, for every $c\in[1,2]$ there exists $p\in[1,2]$ such that $R(\ell^p(2))=c$.
\end{corollary}

\begin{proof}
For $p\in(1,\infty)$ and $\tau\in[0,1]$, let
$$F(p,\tau)=\big(\tau^{p-1}+(1-\tau)^{p-1}\big)^{1/p}\big(\tau^{q-1}+(1-\tau)^{q-1}\big)^{1/q},$$ 
where $q=\frac{p}{p-1}$. Since $p-1$ and $q-1$ are positive numbers, the function $F$ is continuous on $(1,\infty)\times[0,1]$. Hence 
the function $p\mapsto\max_{0\leq\tau\leq1}F(p,\tau)=R(\ell^p(2))$ is continuous as well. Clearly $F(2,\tau)=1$ for all $\tau$.

Fix $\tau\in(0,\frac12)$. If $p\to\infty$, then $q-1\to0$. Hence the second factor of $F(p,\tau)$ tends
to $2$. Since $\big(\tau^{p-1}+(1-\tau)^{p-1}\big)^{1/p}=(1-\tau)^{(p-1)/p}\big(1+(\frac{\tau}{1-\tau})^{p-1}\big)^{1/p}$,
the first factor tends to $1-\tau$. It follows that $\liminf\limits_{p\to\infty}R(\ell^p(2))\geq2(1-\tau)$. Since $\tau$ is
arbitrary and $R(\ell^p(2))\leq2$ by Proposition \ref{prop42}, we get $\lim\limits_{p\to\infty}R(\ell^p(2))=2$. Since
$R(\ell^p(2))=R(\ell^q(2))$, also $\lim\limits_{p\to1}R(\ell^p(2))=2$.

By Example \ref{ex410}, $R(\ell^1(2))=2$. Let $c\in[1,2)$. There is $p_0\in(1,2)$ with $R(\ell^{p_0}(2))>c$, and
$R(\ell^2(2))=1\leq c$. By the intermediate value theorem, $R(\ell^p(2))=c$ for some $p\in[p_0,2]$.
\end{proof}

\subsection*{Acknowledgements}
The first author was supported by the Slovenian Research and Innovation Agency program P2-0268. The second author was supported by FCT/Portugal and the Recovery and Resilience Plan (PRR) through projects:\\ UID/04459/2025 (DOI: 10.54499/UID/04459/2025) and\\ UID/PRR/04459/2025  (DOI: 10.54499/UID/PRR/04459/2025).

\subsection*{AI Declaration}
The authors used Claude (Anthropic, model Claude Opus 5.5) during the development of this work to assist with literature search, 
with developing and checking proofs and examples, and with language editing. All results and arguments were independently verified 
and substantially revised by the authors, who take full responsibility for the content and the final version of the manuscript.


\end{document}